\documentclass[12pt, a4paper]{article}
\usepackage{fullpage}
\usepackage{microtype}
\usepackage[utf8]{inputenc}
\usepackage[english]{babel}
\usepackage[T1]{fontenc}
\usepackage{amsmath,amsthm}
\usepackage{amsfonts}
\usepackage{amssymb}
\usepackage[T1]{fontenc}
\usepackage[shortlabels]{enumitem}
\usepackage{hyperref}
\usepackage{tikz-cd}
\usetikzlibrary{calc, graphs, graphs.standard, trees}
\usepackage{tkz-berge}
\usepackage{caption,graphicx}
\usepackage{hyperref}
\usepackage{xcolor}
\usepackage{chngcntr}
\usepackage{apptools}
\usepackage[title]{appendix}
\AtAppendix{\counterwithin{lemma}{section}}

\usepackage[affil-it]{authblk}
\makeatletter
\def\@maketitle{%
  \newpage
  \null
  \vskip 1em%
  \begin{center}%
  \let \footnote \thanks
    {\Large\bfseries \@title \par}%
    \vskip 0.5em%
    {\normalsize
      \lineskip 0.5em%
      \begin{tabular}[t]{c}%
        \@author
      \end{tabular}\par}%
    {\normalsize \@date}%
  \end{center}%
}

\makeatother
\hypersetup{
    colorlinks,
    linkcolor={red!50!black},
    citecolor={blue!50!black},
    urlcolor={blue!80!black}
}

\newcommand{\N}{\mathbb{N}}

\newtheorem{theorem}{Theorem}
\newtheorem*{theorem*}{Theorem}

\theoremstyle{theorem}
\newtheorem{lemma}{Lemma}
\newtheorem*{lemma*}{Lemma}

\theoremstyle{definition}
\newtheorem{definition}{Definition}

\theoremstyle{definition}
\newtheorem*{definition*}{Definition}

\theoremstyle{definition}
\newtheorem*{example*}{Example}

\theoremstyle{theorem}

\newtheorem*{conjecture*}{Conjecture}

\usepackage{newpxtext, newpxmath}
\title{Equality cases for a bound on the chromatic number}
\author{
Ho Boon Suan, Joel Tan Junyao \& Zhang Xiaorui

}

\date{}

\begin{document}
\maketitle
\abstract{
It is known that the inequality \[\frac{\chi(G)(\chi(G)-1)}{2} + |V| - \chi(G) \leq |E|\] holds for all connected graphs \cite{Soto}, where $\chi(G)$ denotes the chromatic number of $G$. We prove that equality holds whenever the graph consists of a complete graph or an odd cycle, together with finitely many trees attached to its vertices.
} 
\\

Throughout the article, we will only consider {\it simple} graphs (without loops or multiple edges). 
We begin by stating Lemma 1 of \cite{Soto}.
\begin{lemma} \label{lemma:mnbound}
For any connected graph $G = (V, E)$, we have the inequality
\begin{equation} \label{eqn:mnbound}
\frac{\chi(G)(\chi(G)-1)}{2} + |V| - \chi(G) \leq |E|.
\end{equation}
\end{lemma}
We wish to investigate graphs for which (\ref{eqn:mnbound}) is an equality. We begin with a definition.
\begin{definition}
A connected graph that is either
\begin{enumerate}[A.]
\item A complete graph $K_n$ where $n \geq 1$ with finitely many trees attached to each vertex, or
\item An odd cycle $C_{2m+1}$ where $m \geq 1$ with finitely many trees attached to each vertex
\end{enumerate}
is said to be of type A or B respectively.
If a graph $G_A$ or $G_B$ is of type A or B, we will occasionally refer to the complete graph or odd cycle it contains as its underlying graph.

\end{definition}

We have the following result:

\begin{theorem} \label{thm:cool}
A connected graph $G = (V, E)$ satisfies the equality 
\begin{equation} \label{eqn:thm}
\frac{\chi(G)(\chi(G)-1)}{2} + |V| - \chi(G) = |E|
\end{equation}
if and only if it is of type A or B.
\end{theorem}

We will need the following lemma in the proof of \autoref{thm:cool}.

\begin{lemma} \label{lemma:big}
Let $G_n$ be a connected graph with $n$ vertices and $m_n$ edges such that the equality (\ref{eqn:thm}) holds, and suppose that $v \in V(G_n)$ is such that $G_{n-1} := G_n - v$ is connected. Suppose further that \autoref{thm:cool} holds for all $k \leq n$. Then, the following hold:
\begin{enumerate}[(I)]
\item At least one of the following statements are true:
\begin{enumerate}[1.]
\item $d_{G_n}(v) = 1$ and $G_n$ is of type A or B, or
\item $d_{G_n}(v) = \chi(G_{n-1})$ and $G_{n-1}$ is of type A or B;
\end{enumerate}
\item If $d_{G_n}(v) = \chi(G_{n-1}) \leq 2$, then $G_n$ is of type A or B.
\end{enumerate}
\end{lemma}
\begin{proof}
We proceed by induction on the number of vertices. 
There is only one connected graph $G_2$ with two vertices, and it is straightforward to verify that the equality (\ref{eqn:thm}) holds for it. Let $v \in V(G_2)$. Then $G_1 := G_2 - v$ is the graph on one vertex with $\chi(G_1) = 1$. That (I2) and (II) hold is clear, so we are done.

Suppose inductively that the result holds for all $k \leq n$, and let $G_{n+1}$ be a graph with $n+1$ vertices and $m_{n+1}$ edges such that equality in (\ref{eqn:thm}) holds.
Suppose $v \in V(G_{n+1})$ is such that $G_n := G_{n+1} - v$ is connected, and let $m_n$ denote the number of edges of $G_n$. Then, by \autoref{lemma:mnbound}, there exists a non-negative integer $c$ such that
\begin{equation} \label{eqn:n}
\frac{\chi(G_n)(\chi(G_n)-1)}{2} + n - \chi(G_n) + c = m_n.
\end{equation}
Since equality in (\ref{eqn:thm}) holds for $G_{n+1}$ by hypothesis, we have
\begin{equation} \label{eqn:n+1}
\frac{\chi(G_{n+1})(\chi(G_{n+1})-1)}{2} + (n+1) - \chi(G_{n+1}) = m_{n+1}.
\end{equation}

As the addition of a vertex may either increase the chromatic number by one or leave it unchanged, we must consider these cases separately.

{\bf Case 1: $\boldsymbol{\chi(G_{n+1}) = \chi(G_n)}$.} Subtracting (\ref{eqn:n}) from (\ref{eqn:n+1}), we find that \[ d_{G_{n+1}}(v) = m_{n+1} - m_n = 1 - c. \] Since $G_{n+1}$ is connected, we must have $d_{G_{n+1}}(v) > 0$, so that $c = 0$ and $d_{G_{n+1}}(v) = 1$. It follows that equality as in (\ref{eqn:thm}) is satisfied for $G_n$, so that it is of type A or B by hypothesis. Since $d_{G_{n+1}}(v) = 1$, it follows from \autoref{lemma:2} (see Appendix) that $G_{n+1}$ is type A or B. That is, (I1) holds, so that (I) holds. Since $d_{G_{n+1}}(v) = 1$, the antecedent in (II) holds iff $\chi(G_n) = 1$. But that means that $G_n$ is the graph on one vertex, so that $G_{n+1}$ must be the connected graph on two vertices. It follows that $\chi(G_{n+1}) = 2 \neq \chi(G_n)$, a contradiction. Thus (II) holds vacuously. 

{\bf Case 2: $\boldsymbol{\chi(G_{n+1}) = \chi(G_n) + 1}$.}
We claim that
\[d_{G_{n+1}}(v) \geq \chi(G_n).\]
Indeed, suppose for contradiction that $d_{G_{n+1}}(v) < \chi(G_n)$. Then, given a $\chi(G_n)$-coloring of $G_n$, the vertices adjacent to $v$ would not use all $\chi(G_n)$ colors, so there would be an unused color one may assign to $v$. But this gives us a $\chi(G_n)$-coloring of $G_{n+1}$, a contradiction.

By adding $\chi(G_n)$ to both sides of (\ref{eqn:n}), we find that
\[\frac{\chi(G_n)(\chi(G_n)+1)}{2} + n - \chi(G_n) + c = m_n + \chi(G_n).\]
Since $\chi(G_{n+1}) = \chi(G_n) + 1$ by assumption, it follows from (\ref{eqn:n+1}) that
\begin{align*}
m_{n+1} 
&= \frac{\chi(G_{n+1})(\chi(G_{n+1})-1)}{2} + (n+1) - \chi(G_{n+1}) \\
&= \frac{\chi(G_n)(\chi(G_n) + 1)}{2} + n - \chi(G_n) \\
&= m_n + \chi(G_n) - c,
\end{align*}
so that
\[  d_{G_{n+1}}(v) = m_{n+1} - m_n = \chi(G_n) - c.\]
Since $\chi(G_{n+1}) > \chi(G_n)$, it follows that $d_{G_{n+1}}(v) \geq \chi(G_n)$. But since $\chi(G_n) - d_G(v) = c \geq 0$, it follows that $d_{G_{n+1}}(v) = \chi(G_n)$ and that $c = 0$. Thus, equality as in (\ref{eqn:thm}) holds for $G_n$, so that $G_n$ is of type A or B by hypothesis. It follows that (I2) holds, so that (I) holds.

We now prove (II). If $d_{G_{n+1}}(v)=\chi(G_n) \leq 2$, then there are two cases to consider.

{\bf Case 2.1: $\boldsymbol{\chi(G_n) = 1}$.} In this case, $G_n$ must be the graph on one vertex, so that $G_{n+1}$ is the connected graph on two vertices. It follows that $G_{n+1}$ is of type A, so that (II) holds.

{\bf Case 2.2: $\boldsymbol{\chi(G_n) = 2}$.} By (I), it follows that $G_n$ is of type A or B.

\setlength{\leftskip}{15pt}

\indent {\bf Case 2.2.1: $\boldsymbol{G_n}$ is of type A with underlying graph $\boldsymbol{K_m}$ with $\boldsymbol{m \leq 2}$.}  In this case, $G_n$ is a tree. In $G_{n+1}$, we claim that the vertex $v$ must be adjacent to exactly two vertices of $G_n$. Indeed, we have $\chi(G_{n+1}) = \chi(G_n) + 1 = 3$ by hypothesis. If $v$ were adjacent to fewer vertices, then $\chi(G_{n+1}) <3$, and similarly we would have $\chi(G_{n+1}) > 3$ if $v$ were adjacent to more than two vertices.

Since $v$ is attached to two vertices of $G_n$, exactly one cycle is present in $G_{n+1}$. If this cycle has odd length, then $G_{n+1}$ is of type B (with the odd cycle being the underlying graph), so we are done. If this cycle has even length, then $\chi(G_{n+1}) = 2$ by \autoref{lemma:1.5}, so that $\chi(G_{n+1}) \neq \chi(G_n) + 1$, a contradiction.

\indent {\bf Case 2.2.2: $\boldsymbol{G_n}$ is of type A with underlying graph $\boldsymbol{K_m}$ with $\boldsymbol{m \geq 3}$.} By \autoref{lemma:1.5}, we find that $\chi(G_n) = m \neq 2$, a contradiction.

\indent {\bf Case 2.2.3: $\boldsymbol{G_n}$ is of type B.} Similarly, we apply \autoref{lemma:1.5} to show that $\chi(G_n) = 3 \neq 2$, a contradiction.

\setlength{\leftskip}{0pt}
\end{proof}

We are now able to prove \autoref{thm:cool}.

\begin{proof}[Proof of \autoref{thm:cool}]
We proceed by induction on the number of vertices. There is only one connected graph with two vertices, and it is straightforward to verify that the equality (\ref{eqn:thm}) holds for it, as well as that it is of type A (considered as a tree attached to $K_1$). 

Suppose that the result holds for all $k \leq n$ for some $n \in \N$, and let $G_{n+1}$ be a graph with $n+1$ vertices and $m_{n+1}$ edges such that equality in (\ref{eqn:thm}) holds. If $G_{n+1}$ is of type A or B, we are done. If not, by Proposition 1.4.1 of \cite{Diestel}, we may label the vertices $v_1, v_2, \dots, v_{n+1}$ in such a way that the subgraph consisting of vertices $v_1, \dots, v_k$ is connected for $1 \leq k \leq n+1$. It follows that the graph $G_n := G_{n+1} - v_{n+1}$ with $m_n$ edges is a connected graph with $n$ vertices.

By \autoref{lemma:big}, we have at least one of the following is true:

\begin{enumerate}[(I1)]
\item $d_{G_{n+1}}(v_{n+1})=1$ and $G_{n+1}$ is type A or B,
\end{enumerate}
in which case we are done, or
\begin{enumerate}[(I2)]
\item $d_{G_{n+1}}(v_{n+1})=\chi(G_n)$.
\end{enumerate}

In this case, we must consider various cases depending on $\chi(G_n)$. If $\chi(G_n) \leq 2$, then $G_{n+1}$ is of type A or B by (II) of \autoref{lemma:big}, so we are done. It suffices to prove the theorem for $\chi(G_n) \geq 3$.

\indent {\bf Case 1: $\boldsymbol{G_n}$ is of type A with no trees attached.} That is, $G_n = K_n$. In this case, $\chi(G_{n+1}) = \chi(K_n) + 1 = n + 1$, so that $G_{n+1} = K_{n+1}$ (since every $k$-chromatic graph has a $k$-chromatic subgraph of minimum degree at least $k-1$ (\cite{Diestel}, Lemma 5.2.3)). That is, $G_{n+1}$ is of type A.

\indent {\bf Case 2: $\boldsymbol{G_n}$ is of type B with no trees attached.} That is, $G_n = C_n$.

\setlength{\leftskip}{15pt}

\indent {\bf Case 2.1: $\boldsymbol{G_n = C_3}$.} Since $C_3 = K_3$, this is just Case 1.

\indent {\bf Case 2.2: $\boldsymbol{G_n = C_{2m+1}}$ for $\boldsymbol{m \geq 2}$.} In this case, there must exist some $v \in V(G_n)$ that is not adjacent to $v_{n+1}$. Define a coloring such that $v$ and $v_{n+1}$ get mapped to the same color $1$, and color the rest of the vertices in the odd cycle by alternating between colors $2$ and $3$. This shows that $\chi(G_{n+1}) \leq 3$, so that $\chi(G_{n+1}) \neq \chi(G_n) + 1 = 4$, a contradiction. 

\setlength{\leftskip}{0pt}

\indent {\bf Case 3: $\boldsymbol{G_n}$ contains trees.} In this case, there exists some $v' \in V(G_n)$ with $d_{G_n}(v') = 1$. We claim that $G_n - v'$ is connected. Indeed, suppose for contradiction that there exists two vertices $a$, $b \in G_n - v'$ such that there is no path from $a$ to $b$. However, $G_n$ is connected, so there exists a path from $a$ to $b$ in $G_n$. It follows that this this path must pass through $v'$, so that $d_{G_{n}}(v') \geq 2$, a contradiction. 

\setlength{\leftskip}{15pt}

\indent {\bf Case 3.1: $\boldsymbol{d_{G_{n+1}}(v_{n+1}) = 1}$.} If $G_n$ is of type A or B, direct application of \autoref{lemma:2} implies that $G_{n+1}$ is of type A or B. Otherwise, by (I) of \autoref{lemma:big}, we have $G_n$ or $G_{n+1}$ is of type A or B. If $G_n$ is of type A or B, we apply \autoref{lemma:2} again to conclude that $G_{n+1}$ is of type A or B, as desired.

\indent {\bf Case 3.2: $\boldsymbol{d_{G_{n+1}}(v_{n+1})>1}$.} We have that $v_{n+1}$ must be connected to some vertex of $G_n - v'$ in $G_{n+1}$. It follows that $G := G_{n+1} - v' = (G_n - v') \cup \{v_{n+1}\}$ is connected, so that we may apply \autoref{lemma:big} to $G$ to conclude that

\begin{enumerate}[(1)]
\item $G_{n+1}$ is type A or B, in which case we are done; or
\item $d_{G_{n+1}}(v')=\chi(G)$ and $G$ is type A or B.
\end{enumerate}

However, in the second case, we know that $d_{G_{n}}(v')=1$ implies that $d_{G_{n+1}}(v') \leq 2$. By (II) of \autoref{lemma:big}, $G_{n+1}$ is of type A or B. 

\setlength{\leftskip}{0pt}

In all cases, it follows that $G_{n+1}$ is of type A or B, so that the result follows from induction.

We now prove the converse. Let $G$ be a graph of type A. Then either its underlying graph is $K_1$, in which case $G$ is either the graph on one vertex for which the result is quickly established by a straightforward computation, or a tree with $\chi(G) = 2$ for which the result is established by noting that trees always satisfy $|V| - |E| = 1$; or that its underlying graph is $K_ m$ for $m \geq 2$, in which case we can apply \autoref{lemma:1.5} to find that $\chi(G) = m$. Since $|E(K_m)| - |V(K_m)| = m(m-1)/2 - m$, and since every vertex added afterwards also adds exactly one edge, it follows that
\[\frac{\chi(G)(\chi(G)-1)}{2} - \chi(G) = \frac{m(m-1)}{2} - m = |E(G)| - |V(G)|\] 
as desired. On the other hand, given a graph $G$ of type B, its underlying odd cycle $C_{2m+1}$ has $|E(C_{2m+1})| - |V(C_{2m+1})| = 0$. Since the addition of vertices that form trees adds exactly one edge for each added vertex, we find that $|E(G)| - |V(G)| = 0$. By \autoref{lemma:1.5}, $\chi(G) = 3$, so that \[\frac{\chi(G)(\chi(G)-1)}{2}-\chi(G) = \frac{3(3-1)}{2} - 3 = 0 = |E(G)| - |V(G)|.\]
The converse follows.

\end{proof}

\newpage
\begin{appendix}
\begin{appendices}
\section{Technical lemmas}

\begin{lemma} \label{lemma:1.5}
Let $G_A$ and $G_B$ be type A and B graphs respectively, where $K_m$ is the underlying graph for $G_A$. Let $G_C$ be a graph defined similarly to a type B graph, except that its underlying cycle is even. Then the following statements hold:
\begin{enumerate}[(1)]
\item $\chi(G_A) = m$ if $m \geq 2$,
\item $\chi(G_B) = 3$, and
\item $\chi(G_C) = 2$.
\end{enumerate}
\end{lemma}
\begin{proof}
Let $C_{2n+1}$ be the underlying graph of $G_B$. Since odd cycles (with $n \geq 1$) have chromatic number equal to $3$, it suffices to show that the addition of any number of trees to any of the vertices of $C_{2n+1}$ does not change its chromatic number. Indeed, since any tree with more than two vertices has a chromatic number equal to $2$, and since $\chi(G_B) \geq \chi(C_{2n+1})$, the result follows. The proofs for statements (1) and (3) are very similar, and only require noting that $\chi(K_m) = m$ and $\chi(C_{2n}) = 2$.
\end{proof}

\begin{lemma} \label{lemma:2}
Let $G$ be a graph. If $v \in V(G)$ has $d_G(v) = 1$ and if $G - v$ is of type A, then $G$ is of type A. The result also holds if we replace type A with type B.
\end{lemma}
\begin{proof}
Let $G$ be a graph, and let $v \in V(G)$. Suppose that $d_G(v) = 1$ and $G-v$ is of type A. Since $d_G(v) = 1$, it follows that $v$ is either adjacent to a vertex of the underlying graph or to a vertex of a tree. 
In the first case, $v$ forms a tree at that vertex so that $G$ is type A. 
In the second case, $v$ becomes part of a tree at a vertex of the underlying graph, so that $G$ is type A. 
The proof works exactly the same for type B graphs.
\end{proof}

\end{appendices}
\end{appendix}


\begin{thebibliography}{9}
\bibitem{Soto} M. Soto, A. Rossi and M. Sevaux. Three new upper bounds on the chromatic number. \textit{Discrete Applied Mathematics}, 159 (2011), pp. 2281--2289.

\bibitem{Diestel} R. Diestel. Graph Theory, 5th Ed. (2017). {\it Graduate Texts in Mathematics} 173, Springer-Verlag.
\end{thebibliography}
\end{document}